\documentstyle[12pt]{article}
\begin{document}
\title{Morphisms of Multiplicative Unitaries}
\author{Chi-Keung Ng\footnotemark}
\addtocounter{footnote}{1}\footnotetext{This work is partially 
supported by the Croucher Foundation.
\\
Mathematics Subject Classification 46L55.}
\date{}
\maketitle
\begin{abstract}
In this paper, we will give a natural definition for morphisms 
between multiplicative unitaries. 
We will then discuss 
some equivalences of this definition and some interesting 
properties of them. 
Moreover, we will define normal sub-multiplicative 
unitaries for multiplicative unitaries of discrete type and prove an 
imprimitivity type theorem for discrete multiplicative unitaries. 
\end{abstract}
\par\medskip
\noindent {\bf 0. Introduction}
\par
\medskip
In $[2]$, Baaj and Skandalis defined multiplicative unitaries and 
showed that they are nice generalisation of locally compact groups. They 
also showed that the Woronowicz $C^{*}$-algebras (which can also be 
considered as compact quantum groups) can be included in the 
consideration of multiplicative unitaries of compact type. 
However, there isn't any notion of morphism being defined so far. 
In [8], Wang define morphisms between compact quantum groups as Hopf 
*-homomorphisms between the underlying Woronowicz $C^{*}$-algebras. 
However, for a given 
multiplicative unitary $V$, we can associate with it 
four Hopf $C^{*}$-algebras 
(if the multiplicative unitary is good), 
namely, $S_{V}, \hat{S}_{V}, (S_{V})_{p}$ 
and $(\hat{S}_{V})_{p}$. It is not clear which of 
the Hopf *-homomorphisms 
between these Hopf $C^{*}$-algebras should be used as a candidate for 
the morphisms.

In this paper, we will investigate a natural notion called 
``birepresentation'' and show that it is a good candidate for the morphisms 
between multiplicative unitaries. More precisely, given two 
multiplicative 
unitaries U and V, we define morphisms from V to U to be the collection of 
all U-V-birepresentations. Some of the Hopf *-homomorphisms between the 
Hopf $C^{*}$-algebras defined by $U$ and $V$ 
are equivalent to $U$-$V$-birepresentations (see Theorem 4.9). 
We also find another equivalence that birepresentations are 
in one to one correspondence with the ``mutual coactions'' 
(see Definition 3.13) 
between those Hopf $C^{*}$-algebras. 
We also investigate the crossed products of the 
different coactions arise from a ``morphism'' and 
show that $(S_{V})_{p}\times _{\delta ,\max }(\hat{S}_{V})_{p} \cong 
(\hat{S}_{V})_{p}\times _{\hat{\delta }^{op},\max }(S_{V})_{p}^{op}$.

Now we obtain a category of multiplicative unitaries 
(those satisfy some good property).
It contains the locally compact groups as a full subcategory. 
However, it seem not easy to define kernels of these morphisms. 
By looking at the the case of discrete groups, 
we can define kernels of 
morphisms between multiplicative unitaries of discrete type (which is 
really a kernel in the categorical sense).

Finally, from the definition of kernel, we can define normal 
sub-multiplicative unitaries of multiplicative unitaries of discrete type. 
We then prove an imprimitivity type theorem for this setting.

We would like to thank Prof. Zhong-Jin Ruan for indicating 
a mistake in the original version of this paper.
\par
\medskip
\par
\medskip
\noindent {\bf 1. Preliminary and Notations}
\par
\medskip
The notations in this paper mainly follow from those of [2] and [5]. We also 
assume basic definitions and results from these two papers. 
\par
\medskip
\noindent {\bf Definition} 1.1{\bf :} Let $(A, \delta )$ be a Hopf 
$C^{*}$-algebra. 

\noindent (a) A *-subalgebra $B$ of
$M(A)$ is called a {\it Hopf $C^{*}$-subalgebra} of $A$ if
\par
\noindent (i) there exists an approximate unit $\{e_{i}\}$ of $B$ such that 
$e_{i}$ converges strictly
to 1 in $M(A)$;
\par
\noindent (ii) $\delta (B) \subseteq M(B\otimes B);$
\par
\noindent (iii) the restriction, $\epsilon $, of $\delta $ in $B$ is a 
comultiplication on $B$.
\par
\noindent (b) Let $A^{op}$ be the $C^*$-algebra $A$ with a comultiplication 
$\delta^{op}$ defined by $\delta^{op} = \sigma\circ\delta$ (where $\sigma$ 
is the flip of variables).
\par
\medskip
Note that condition (iii) in 1.1(a) means that $\epsilon (B)(1\otimes B) 
\subseteq B\otimes B$ and condition 
(ii) makes sense because of (i). It is easy to see that $\delta^{op}$ 
is a comultiplication on $A$ and $(A^{op})^* = 
(A^*)^{op}$ (algebra with an opposite multiplication). We recall 
that a {\it Hopf *-homomorphism} $\phi $ from
a Hopf $C^{*}$-algebra $(A, \delta )$ to another Hopf $C^{*}$-algebra $(B, 
\epsilon )$ is a non-degenerate
*-homomorphism from $A$ to $M(B)$ such that $(\phi \otimes \phi )\circ 
\delta = \epsilon \circ \phi $.
\par
\medskip
\noindent {\bf Lemma} 1.2{\bf :} Let $(A, \delta )$ and $(B, \epsilon )$ be 
two Hopf $C^{*}$-algebras and $\phi $ be a Hopf 
*-homomorphism from $A$ to $M(B)$. Then $B_{0} = \phi (A)$ is a Hopf 
$C^{*}$-subalgebra of $B$.
\par
\medskip
\noindent {\bf Lemma} 1.3{\bf :} Let $A$ and $B$ be $C^{*}$-algebras. If 
$\phi $ and $\psi $ are non-degenerate
*-homomorphisms from $A$ to $M(B)$ and from $B$ to $M(A)$ respectively such 
that $\phi \circ \psi $
and $\psi \circ \phi $ are identity maps, then $\phi $ is an isomorphism 
from $A$ to $B$.
\par
\noindent {\bf Proof:} We first note that $\phi (A)$ is an ideal of $M(B)$ 
since $\psi (\phi (a)m) = a\psi (m) \in A$
for any $a\in A$ and $m\in M(B)$. Let $\{e_{i}\}$ be an approximate unit of 
A. Then for any
$b\in B$, $\phi (e_{i})b \in \phi (A)$ will converge to $b$ and hence $B 
\subseteq \phi (A)$.
\par
\medskip
\noindent {\bf Definition} 1.4{\bf :} Let $(R, \epsilon )$ and 
$(S, \delta )$ be Hopf $C^{*}$-algebras. An unitary $U\in M(R\otimes S)$
is said to be a {\it unitary $R$-$S$-birepresentation} if 
$(id\otimes \delta )(U)=U_{12}U_{13}$ and $(\epsilon \otimes id)(U) =
U_{13}U_{23}$.
\par
\medskip
\noindent {\bf Lemma} 1.5{\bf :} Let $(S, \delta )$ and $(T,\epsilon )$ be 
Hopf $C^{*}$-algebras. Let $w$ and $v$ be unitary
co-representations of $S$ and $T$ respectively on the same Hilbert 
space $H$ and let $u = w^{\sigma } \in 
M(S\otimes {\cal K}(H))$ (where $\sigma $ means the flip of the two 
variables). If $X$ is a unitary in
$M(S\otimes T)$ such that $u_{12}X_{13}v_{23} = v_{23}u_{12}$, then $X$ is a 
unitary $S$-$T$-birepresentation.
\par
\noindent {\bf Proof:} Applying $(id\otimes id\otimes \epsilon )$ on the 
equation, we have that $u_{12}(id\otimes id\otimes \epsilon )
(X_{13})v_{23}v_{24} =
v_{23}v_{24}u_{12} = v_{23}u_{12}X_{14}v_{24}$. Thus, 
$(id\otimes id\otimes \epsilon )(X_{13}) = $
$u^{*}_{12}v_{23}u_{12}X_{14}v^{*}_{23} = $
$X_{13}v_{23}X_{14}v^{*}_{23} = $
$X_{13}X_{14}$. Similarly, $(\delta \otimes id)(X) = X_{13}X_{23}$.
\par
\medskip
\noindent {\bf Definition} 1.6{\bf :} Let $X$ be a unitary 
$S$-$T$-birepresentation. Let $w$ and $v$ be
unitary co-representations of $S$ and $T$ respectively on the same 
Hilbert space and let 
$u = w^{\sigma }$. Then $(u, v)$
is said to be a {\it covariant pair for $X$} if $u$ and $v$ satisfy the 
condition in
the previous lemma.
\par
\medskip
We now recall the following definitions from [1].
\par\medskip
\noindent {\bf Definition} 1.7{\bf :} Let $V$ be a multiplicative unitary 
on a Hilbert space $H$. Then 
\par
\noindent (a) $V$ is said to be semi-regular if the 
norm closure of the set 
$\{ ($id$\otimes\omega)(\Sigma V): \omega\in {\cal L}(H)_* \}$ contains 
the set of all compact operators ${\cal K}(H)$. Moreover, $V$ is said to 
be semi-biregular if it is regular and the norm closure of the set 
$\{ (\omega\otimes$id$)(\Sigma V): \omega\in {\cal L}(H)_* \}$ contains 
${\cal K}(H)$ as well. 
\par
\noindent (b) $V$ is said to be balanced if there exists a unitary 
$U\in {\cal L}(H)$ such that 
\par
\noindent (i) $U^2=I_H$;
\par
\noindent (ii) the unitary $\hat{V} = 
\Sigma(U\otimes 1)V(U\otimes 1)\Sigma$ is multiplicative. 
\par\medskip
\noindent {\bf Remark} 1.8{\bf :}
For simplicity, we will call a multiplicative unitary 
{\it semi-irreducible} if it is both semi-regular and balanced.
\par
\medskip
\noindent {\bf Proposition} 1.9{\bf :} Let $S$ and $T$ be Hopf 
$C^{*}$-algebras and $\phi $ be a Hopf
*-homomorphism from $S$ to $M(T)$. If $\epsilon $ is a coaction on a $
C^{*}$-algebra $A$ by $S$, then
$\delta = (id\otimes \phi )\circ \epsilon $ is a coaction on $A$ by $T$.
\par
\noindent {\bf Proof:} The coaction identity follows 
easily from the fact that $\phi $
respects the comultiplications. It remains to show that 
$\delta (A) \subseteq 
\tilde{M} (A\otimes T)$. Let $(u_{i})$ be an approximate unit of S. For any 
$a\in A$ and $t\in T$, $t_{i}=\phi (u_{i})t$
converges to $t$ in norm and so $\delta (a)(1\otimes t_{i})$ converges to 
$\delta (a)(1\otimes t)$ in norm. Now,
$\delta (a)\cdot (1\otimes t_{i}) = (id\otimes \phi )[\epsilon (a)(1\otimes 
u_{i})]\cdot (1\otimes t) \in A\otimes T$.
\par
\medskip
\noindent {\bf Proposition} 1.10{\bf :} Let $A, S, T, \phi , \epsilon $ and 
$\delta $ be the same as in Proposition 1.9.
Suppose that the crossed products $A\times _{\epsilon ,\max }\hat{S}$ and 
$A\times _{\delta ,\max }\hat{T}$ exist. Then there
exists a *-homomorphism $\Phi $ from $A\times _{\epsilon ,\max }\hat{T}$ 
to $M(A\times _{\delta ,\max }\hat{S})$.
\par
\noindent {\bf Proof:} Let $(B, \psi , u)$ be a covariant pair for $(A, S, 
\epsilon )$ and let $v = (id\otimes \phi )(u)$.
Then $(B, \psi , v)$ is a covariant pair for $(A, T, \delta )$ and the 
proposition follows
from the definition of crossed product (see $[5, 2.11(b)])$.
\par\medskip
\par\medskip
\noindent {\bf 2. Basic Mulitplicative Unitaries}
\par
\medskip
The aim of this section is to find some basic assumptions on the
multiplicative unities such that the results in this paper holds. We
will show that the semi-irreducible multiplicative unitaries 
and regular mulitplicative unitaries both satisfy these basic
assumptions (the manageable multiplicative unitaries ``almost''
satisfy these assumptions, at least when they are either amenable or 
co-amenable). 
\par\medskip
\noindent {\bf Definition} 2.1: Let $V$ be a multiplicative unitary. 
Then $V$ is a called a {\it $C^*$-multiplicative unitary} if for any 
representation $X$ and corepresentation $Y$ of $V$ on $K$
and $L$ respectively, 
\par
\noindent (i) $\hat S_X = \overline 
{\{ (id\otimes \omega)(X): \omega\in {\cal L}(H)_* \}}$ and $S_Y=  
\overline {\{ (\omega\otimes id)(Y): \omega\in {\cal L}(H)_* \}}$ 
are both $C^*$-algebras;
\par
\noindent (ii) $X\in M(\hat S_X\otimes S_V)$ and 
$Y\in M(\hat{S}_V\otimes S_Y)$.
\par\medskip
Basic examples of $C^*$-multiplicative unitaries are regular
multiplicative unitaries, semi-irreducible multiplicative
unitaries (see Remark 1.8) and manageable multiplicative unitaries. 
\par\medskip
\noindent {\bf Remark} 2.2: 
(a) By the arguement in [11, section 5], if $V$ is a
$C^*$-multiplicative unitary, then 
$S_V$ and $\hat S_V$ are both Hopf
$C^*$-algebras with coaction $\delta$ and $\hat \delta$ respectively.
Moreover, $(id\otimes \delta)(X)=X_{12}X_{13}$ and $(\hat \delta
\otimes id)(Y)=Y_{13}Y_{23}$ if $X$ and $Y$ are representation and 
corepresentation of $V$ respectively. Furthermore, by using the 
arguement in that section, we can also show that the closure of 
$\{ (\omega\otimes id\otimes id)(Y_{12}Y_{13})(1\otimes s): 
\omega\in {\cal L}(H)_*, s\in S_Y \} = S_Y\otimes S_Y$ for any 
corepresentation $Y$.
\par
\noindent (b) If $V$ is a $C^*$-multiplicative unitary, then 
by the same arguement as in $[2, A6]$, $[2, A6(a)-(d)]$ hold for $V$. 
Moreover by part (a) above, $[2, A6(e)]$ holds as well.
\par\noindent 
(c) We also note that all the main results in [5] hold for
$C^*$-multiplicative unitaries (actually, except [5, 3.7 \& 
3.15] which involve the Takesaki-Takai type duality).
\par\medskip
Let $V'$ and $V''$ be as defined in [2, A6]. We call $V'$ and $V''$ the
universal representation and the universal corepresentation of $V$
respectively. We also need the following technical assumption. 
\par\medskip
\noindent {\bf Definition} 2.3: Let $V$ be a $C^*$-multiplicative unitary.
\par
\noindent (a) Let $V_p$ be a unitary in $M((\hat S_V)_p\otimes
(S_V)_p)$. Then $V_p$ is said to be a {\it universal birepresentation}
of $V$ if $V'_{12}(V_p)_{13}V''_{23}=V''_{23}V'_{12}$ in 
$M((\hat S_V)_p\otimes {\cal K}(H)\otimes (S_V)_p)$.
\par
\noindent (b) A $C^*$-mulitplicative unitary $V$ is said to be 
{\it basic} if there exists a universal birepresentation for $V$.
\par\medskip
\noindent {\bf Remark} 2.4: (a) It is clear that if a basic
multiplicative unitary $V$ is amenable (respectively, coamenable), 
then $V_p=V''$ (respectively, $V_p=V'$) exists 
and $V$ is basic.
\par
\noindent (b) By Lemma 1.5, if the universal
birepresentation of $V$ exists, it is a unitary 
$(\hat S_V)_p$-$(S_V)_p$-birepresentation. Moreover, it is clear that 
$(id\otimes L_V)(V_p)=V'$ and $(\rho_V\otimes id)(V_p)=V''$. 
\par\medskip
We are going to show that $V_p$ exists in good case. 
Note that we can also deduce the existence of $V_p$ from $[2, A8]$ 
but irreducibility is required there. 
We first recall the set ${\cal C}(V)= \{ (id\otimes \omega)
(\Sigma V): \omega\in {\cal L}(H)_* \}$ from [2]. 
Note that the idea of the proof of the following lemma is from [2,3.6(c)].
\par\medskip
\noindent {\bf Lemma} 2.5: Let $V\in {\cal L}(H\otimes H)$ be a 
multiplicative unitary and $X$ and $Y$ are a representation 
and a corepresention of $V$ on $K$ and $L$ respectively. Let
$W=X^*_{12}Y_{23}X_{12}Y^*_{23} \in {\cal L}(K\otimes H\otimes L)$.
Then $(1\otimes c\otimes 1)W=W(1\otimes c\otimes 1)$ for any $c\in
{\cal C}(V)$. Consequently, if $\overline{{\cal C}(V)}^{weak}=
{\cal L}(H)$, then $W$ is of the form $W=Z_{13}$ for some unitary
$Z\in {\cal L}(K\otimes L)$. 
\par
\noindent {\bf proof:} We first note that 
$X^*_{12}Y_{24}X_{12}Y_{24}^*\Sigma_{23}V_{23}=$ 
$\Sigma_{23}X^*_{13}Y_{34}X_{13}Y_{34}^*V_{23}=$
$\Sigma_{23}X^*_{13}Y_{34}X_{13}V_{23}Y_{34}^*Y_{24}^*=$
$\Sigma_{23}X^*_{13}Y_{34}X^*_{12}V_{23}X_{12}Y_{34}^*Y_{24}^*=$
$\Sigma_{23}(X_{12}X_{13})^*V_{23}Y_{24}Y_{34}X_{12}Y_{34}^*Y_{24}^*=$
$\Sigma_{23}V_{23}X^*_{12}Y_{24}X_{12}Y_{24}^*$. 
Now let $c=(id\otimes \omega)(\Sigma V)$. Then 
$W(1\otimes c\otimes 1)=(id\otimes id\otimes \omega\otimes id)
(X^*_{12}Y_{24}X_{12}Y_{24}^*\Sigma_{23}V_{23})=
(id\otimes id\otimes \omega\otimes id)
(\Sigma_{23}V_{23}X^*_{12}Y_{24}X_{12}Y_{24}^*)=
(1\otimes c\otimes1)W$.
The finally part of the proposition is clear.
\par\medskip
\noindent {\bf Proposition} 2.6{\bf :} If $V$ is a $C^*$-multiplicative
unitary such that $\overline{{\cal C}(V)}^{weak}=
{\cal L}(H)$, then $V$ is basic.
\par
\noindent {\bf Proof:} This proposition is clear by putting $X=V'$ and
$Y=V''$ into Lemma 2.5.
\par\medskip
\noindent {\bf Corollary} 2.7{\bf :} Every semi-irreducible (respectively, 
regular) multiplicative unitary is basic.
\par\medskip
\noindent {\bf Remark} 2.8{\bf :} It is natural to ask whether we can use 
similar argument as in Lemma 2.5 to prove that manageable multiplicative 
unitaries are basic as well. However, we encounter a difficulty in doing 
so. The difficulty come from the unboundedness of $Q$. It is not hard to 
show that if $(V, Q, \tilde V)$ is a manageable multiplicative unitary 
such that $(f\otimes id\otimes g)(W)(Dom Q)\subseteq Dom Q$ (where $W=
V^{'*}_{12}V^{''}_{23}V^{'}_{12}Y^{''*}_{23}$) for any 
$f\in (\hat S_V)^*_p$ and $g\in (S_V)^*_p$, then $V$ is basic. 
Note that if $Q$ is bounded, then $V$ is regular and $\kappa_V$ is bounded.
\par\medskip
In the remainder of this section, we assume the multiplicative unitary
$V\in {\cal L}(H\otimes H)$ to be basic. 
Let $V^{\top }=\Sigma V^{*}\Sigma 
\in {\cal L}(H\otimes H)$ and 
$V_{p}\in M[(\hat{S} _{V})_{p}\otimes (S_{V})_{p}]$ be the universal 
birepresentation of $V$.
\par
\medskip
\noindent {\bf Lemma} 2.9{\bf :} (i) $S_{V} = \hat{S}_{V^{\top }}$ and 
$\hat{S}_{V} = S_{V^{\top }}$ as *-subalgebras of ${\cal L}(H)$ (in fact, 
$S_{V}^{op} = \hat{S}_{V^{\top }}$ and $ \hat{S}_{V}^{op} = S_{V^{\top }}$
as Hopf $C^*$-algebras);
\par
\noindent (ii) If $W$ is a representation (respectively, co-representation) 
of $V$, then $W^{\top }=
\Sigma W^{*}\Sigma $ is a co-representation (respectively, representation) 
of $V^{\top }$;
\par
\noindent (iii) $(S_{V})_{p}^{op} \cong (\hat{S}_{V^{\top }})_{p}$ and 
$(\hat{S}_{V})_{p}^{op} \cong (S_{V^{\top }})_{p}$ 
(as Hopf $C^{*}$-algebras); 
\par
\noindent (vi) $(V^{\top })^\prime = (V^{\prime\prime *})^{\sigma} $, 
$(V^{\top })^{\prime\prime} = (V^{\prime *})^{\sigma}$ and 
$V^{\top }_{p} = (V^{*}_{p})^{\sigma }$ (where $\sigma$ means the flip
of variables). 
\par
\medskip
We recall the antipode $\kappa _{V}$ from $A_V$ to $S_{V}$ defined by 
$\kappa _{V}((\omega \otimes id)(V)) = (\omega \otimes id)(V^{*})$. 
In the same way, we define the antipode $j_V$ from $(A_{V})_p = 
\{(\omega \otimes id)(V''): \omega \in {\cal L}(H)_*\}$ to $(S_V)_p$ by 
$j_V((\omega\otimes id)(V^{\prime\prime})) = (\omega \otimes id) 
(V^{\prime\prime *})$. Note that $j_V$ is well defined since 
$(\omega \otimes id)(V'') = 0$ implies $\omega = 0$ on $\hat {S}$ and hence 
$\omega^* = 0$ on $M(\hat {S})$ (which implies that 
$\omega^*((id \otimes f)
(V'')) = 0$ for all $f \in (S_V)_p^*$). 
Moreover, we can extend $\kappa _V$ and $j_V$ as follow.
\par
\medskip
\noindent {\bf Lemma} 2.10{\bf :} 
$\kappa_V$ and $j_V$ can be extended to $\tilde {A_{V}} = 
\{(f \otimes id)(V'): f\in (\hat{S}_V)_p^*\}$ and $\widetilde {(A_{V})_p} = 
\{(f \otimes id)(V_p): f\in (\hat{S}_V)_p^*\}$ respectively.
\par
\noindent {\bf Proof:} 
Since the map that send $f\in (\hat{S}_V)_p^*$ to $(f\otimes id)(V')$ is 
injective (see [5, A6]), 
the map $\kappa_V$ that send $(f\otimes id)(V')$ 
to $(f\otimes id)(V'^*)$ is well defined and is clearly an extension of the 
$\kappa_V$ above. Similarly, since the map that send 
$f\in (\hat{S}_V)_p^*$ to $(f\otimes id)(V_p)$ is injective, the extension 
of $j_V$ is also well defined. 
\par
\medskip
\noindent {\bf Proposition} 2.11{\bf :} There is a one to one 
correspondence between unitary
co-representations of $S_{V}$ and those of $(S_{V})_{p}$.
\par
\noindent {\bf Proof:} If $w$ is a unitary co-representation of $S_{V}$, 
then $w_{p}=(\rho _{w}\otimes id)(V_{p})$ is a
unitary co-representation of $(S_{V})_{p}$ (by the Remark 2.4(b)). 
On the other hand, if
$u$ is a unitary co-representation of $(S_{V})_{p}$, then $u_{0}=(id\otimes 
L_{V})(u)$ is a unitary
co-representation of $S_{V}$. Moreover, it is clear that $(id\otimes 
L_{V})(\rho _{w}\otimes id)(V_{p}) =
(\rho _{w}\otimes id)(V^\prime ) = w$. It remains to show that if 
$(id\otimes L_{V})(u_{1}) = (id\otimes L_{V})(u_{2})$, then
$u_{1}=u_{2}$. It follows from exactly the same argument as in $[5, 2.7]$.
\par
\medskip
\noindent {\bf Corollary} 2.12{\bf :} If $\epsilon ^\prime $ is a 
coaction of a $C^{*}$-algebra $A$ by the Hopf $C^{*}$-algebra
$(S_{V})_{p}$, then the full crossed product $A\times _{\epsilon ^\prime 
,\max }(\hat{S}_{V})_{p}$ exists and is a quotient of 
$A\times _{\epsilon ,\max }\hat{S}_{V}$ (where $\epsilon $ is the reduced 
coaction that corresponds to $\epsilon ^\prime $ as defined in the
paragraph before $[5, 2.14])$.
\par
\noindent {\bf Proof:} Using Proposition 2.11 and the same argument 
as in [5, 2.12(a)], we can reformulate the full crossed product of 
$(A, \epsilon ')$ as in [5, 2.12(c)]. Now by a similar argument as in 
[5, 2.13], the full crossed product exists. Since any covariant 
representation of $(A, \epsilon ')$ is a covariant representation 
of $(A, \epsilon)$, it is clear that $A\times _{\epsilon ^\prime 
,\max }(\hat{S}_{V})_{p}$ is a quotient of 
$A\times _{\epsilon ,\max }\hat{S}_{V}$. 
\par
\medskip
\noindent {\bf Proposition} 2.13{\bf :} Let $A$ be a $C^{*}$-algebra and 
$\epsilon$ a coaction on $A$ by $S_{V}$. Let $(B,
\phi , \mu )$ be the full crossed product. Then there is a 
{\it dual coaction} $\bar{\epsilon }$ on $B$ by
$(\hat{S}_{V})_{p}$ such the $\mu $ is equivariant.
\par
\noindent {\bf Proof:} Let $v = (\mu \otimes id)(V^\prime )$. Then $(\phi 
\otimes id)\epsilon (a)\cdot v = v\cdot (\phi (a)\otimes 1)$ for any 
$a\in A$. Now
define the *-homomorphisms $\psi =(id\otimes 1)\circ \phi $ and $\nu =(\mu 
\otimes id)\circ \hat{\delta }_{V}$ from $A$ and $(\hat{S}_{V})_{p}$ 
respectively to $M(B\otimes (\hat{S}_{V})_{p})$. We first show 
that $(\psi ,\nu )$ is a covariant pair for $(A, S_{V}, \epsilon )$. 
In fact, for any $a\in A,
(\psi \otimes id)\epsilon (a) = (\phi \otimes id)\epsilon (a)_{13}$ and 
$(\nu \otimes id)(V^\prime ) = (\mu \otimes id\otimes id)(\hat{\delta }
_{V}\otimes id)(V^\prime ) = v_{13}V_{23}^\prime $.
Hence $(\psi \otimes id)\epsilon (a)\cdot (\nu \otimes id)(V^\prime ) = 
[(\phi \otimes id)\epsilon (a)\cdot v]_{13}V_{23}^\prime = v_{13}\cdot (\phi 
(a)\otimes 1\otimes 1)\cdot V_{23}^\prime =
(\nu \otimes id)(V^\prime )\cdot (\psi (a)\otimes 1)$. Thus, we have a map 
$\bar{\epsilon }$ from $B$ to $M(B\otimes (\hat{S}_{V})_{p})$ such that 
$\psi =
\bar{\epsilon }\circ \phi $ and $\nu = \bar{\epsilon }\circ \mu $. Now for
 any $a\in A$ and $s,t\in (\hat{S}_{V})_{p}$, $\bar{\epsilon }
[\phi (a)\mu (s)]\cdot (1\otimes t) =
(\phi (a)\otimes 1)\cdot (\mu \otimes id)[\hat{\delta }_{V}(s)\cdot 
(1\otimes t)] \in B\otimes (\hat{S}_{V})_{p}$ (as 
$\hat{\delta }_{V}(s)\cdot (1\otimes t) \in 
(\hat{S}_{V})_{p}\otimes (\hat{S}_{V})_{p}$ and $\phi (a)\mu (u)\in B$ for 
all $u\in (\hat{S_{V}})_{p}$). Since $\{\phi (a)\mu (s): a\in A, 
s\in (\hat{S}_{V})_{p}\}$ generates $B$ (by 
$[5, 2.12(b)(3)]$), we have $\bar{\epsilon }(B) \subseteq 
\tilde{M}(B\otimes (\hat{S}_{V})_{p})$. 
It remains to show the coaction identity.
For $a\in A$ and $s\in (\hat{S}_{V})_{p}$, $(\bar{\epsilon }\otimes 
id)\bar{\epsilon }[\phi (a)\mu (s)] = (\bar{\epsilon }\otimes id)[(\phi 
(a)\otimes 1)\cdot (\mu \otimes id)\hat{\delta }_{V}(s)] =
(\phi (a)\otimes 1\otimes 1)\cdot [(\mu \otimes id)\circ 
\hat{\delta }_{V}\otimes 
id]\hat{\delta }_{V}(s) = (\phi (a)\otimes 1\otimes 1)\cdot (\mu \otimes 
id\otimes id)(id\otimes \hat{\delta }_{V})\hat{\delta }_{V}(s)$. On the
other, we have $(id\otimes \hat{\delta }_{V})\bar{\epsilon }[\phi (a)
\mu (s)] = 
(id\otimes \hat{\delta }_{V})[(\phi (a)\otimes 1)\cdot (\mu \otimes 
id)\hat{\delta }_{V}(s)] =
(\phi (a)\otimes 1\otimes 1)\cdot (\mu \otimes id\otimes id)(id\otimes 
\hat{\delta }_{V})\hat{\delta }_{V}(s)$. 
Finally, $\mu $
is equivariant by the definition of $\bar{\epsilon }$.
\par
\medskip
\par
\medskip
\noindent {\bf 3. $U$-$V$-birepresentations}
\par
\medskip
In this section we mainly deal with $C^*$-multiplicative
unitaries. We will discuss basic multiplicative unitaries in the next 
section. We will define 
and study the birepresentation of
two $C^*$-multiplicative unitaries. 
Let $(K, U)$ and $(H, V)$ be $C^*$-multiplicative unitaries and 
$X\in {\cal L}(K\otimes H)$ be a unitary operator. 
\par
\medskip
\noindent {\bf Definition} 3.1{\bf :} $X$ is said to be a 
{\it $U$-$V$-birepresentation} if $X$ is a representation
of $V$ as well as a corepresentation of $U$.
\par
\medskip
Let $X$ be a $U$-$V$-birepresentation. Then there are *-representations 
$L_{X}$ and $\rho _{X}$ of $(S_{U})_{p}$ and $(\hat{S}_{V})_{p}$ 
on $H$ and $K$ respectively. Moreover, we have:
\par
\medskip
\noindent {\bf Proposition} 3.2{\bf :} $L_{X}$ is a Hopf *-homomorphism 
from $(S_{U})_{p}$ to $M(S_{V})$.
Consequently, $S_{X} = L_{X}(S_{U})_{p}$ is a Hopf $C^{*}$-subalgebra of 
$S_{V}$. Moreover, $L_{X}$ preserves antipodes.
\par
\noindent {\bf Proof:} The first statement is clear from the fact that 
$(\rho _{X}\otimes id)(V^{\prime }) = X$ and $V^\prime \in 
M((\hat{S}_{V})_{p}\otimes S_{V})$. The
second one follows from Lemma 1.2. Finally, we need to show that 
$L_{X}\circ j_{U} = 
\kappa _{V}\circ L_{X}$. Note that $L_{X}\circ j_{U}((\omega \otimes 
id)(U^{\prime\prime})) = (\omega \otimes id)(X^{*})$.
Now since $(\rho _{X}\otimes id)(V^\prime ) = X$, $L_{X}((\omega \otimes id)
(U^{\prime\prime})) = (\omega \circ \rho _{X}\otimes id)(V^\prime) \in 
\tilde {A_V}$ (see Lemma 2.10). Moreover, $\kappa _{V}\circ 
L_{X}((\omega \otimes id)(U^{\prime\prime})) = 
(\omega \circ \rho _{X}\otimes 
id)(V^{\prime *}) = (\omega \otimes id)(X^{*})$.
\par
\medskip
\noindent {\bf Remark} 3.3{\bf :} 
(a) Similar things hold for $\rho _{X}$ 
and $\hat{S}_{X}$.
\par
\noindent (b) It is clear that $X$ is a birepresentation if and only if 
$X$ is a unitary $\hat{S}_{U}$-$S_{V}$-birepresentation in 
$M(\hat{S}_{U}\otimes S_{V})$ (see Definition 1.4).
\par
\noindent (c) Let $X$ be a $U$-$V$-birepresentation. Then $X^{\top 
}=\Sigma X^{*}\Sigma $ is a $V^{\top }$-$U^{\top }$-birepresentation.
Moreover, $L_{X^{\top }} = \rho _{X}$ and $\rho _{X^{\top }} =L_{X}$.
\par
\noindent (d) If we borrow Proposition 3.5 below, we have the following: 
{\it Any Hopf-*-homomorphism from $(S_{U})_{p}$ to $M(S_{V})$ 
preserves antipodes}.
\par
\medskip
We are going to give a converse to proposition 3.2. Let us first investigate 
under what condition a corepresentation will be a birepresentation.
\par
\medskip
\noindent {\bf Lemma} 3.4{\bf :} Let $X$ be a co-representation of $U$ on H. 
If $L_{X}$ is Hopf
*-homomorphism from $(S_{U})_{p}$ to $M(S_{V})$, then $X$ is a 
$U$-$V$-birepresentation.
\par
\noindent {\bf Proof:} It is required to show that $X$ is a representation 
of V. For any $\omega \in {\cal L}(K)_{*}$, 
we have $(L_{X}\otimes L_{X})(\delta _{U}(\omega \otimes 
id)(U^{\prime\prime})) = (\omega \otimes id\otimes id)(X_{12}X_{13})$ 
and $\delta _{V}(L_{X}(\omega \otimes id)(U^{\prime\prime})) =
V[(\omega \otimes id)(X)\otimes 1]V^{*} = (\omega 
\otimes id\otimes id)(V_{23}X_{12}V^{*}_{23})$. Since 
${\cal L}(K)_{*}$ separates points of ${\cal L}(K)$, $X_{12}X_{13}V_{23} 
= V_{23}X_{12}$.
\par
\medskip
\noindent {\bf Proposition} 3.5{\bf :} Any Hopf *-homomorphism $\pi $ from 
$(S_{U})_{p}$ to $M(S_{V})$ induces a unique 
$U$-$V$-birepresentation $X$ such that $\pi = L_{X}$. 
Similarly, any Hopf *-homomorphism from $(\hat{S}
_{V})_{p}$ to $M(\hat{S}_{U})$ also induces a $U$-$V$-birepresentation.
\par
\medskip
\noindent {\bf Theorem} 3.6{\bf :} There are one to one correspondences 
between the followings:
\par
\noindent (a) $U$-$V$-birepresentations;
\par
\noindent (b) Hopf *-homomorphisms from $(S_{U})_{p}$ to $M(S_{V})$;
\par
\noindent (c) Hopf *-homomorphisms from $(\hat{S}_{V})_{p}$ to 
$M(\hat{S}_{U})$.
\par\medskip
\noindent {\bf Corollary} 3.7{\bf :} Let $\phi$ be a non-degenerate
representation of $(S_U)_p$ on ${\cal L}(H)$. Then $\phi((S_U)_p)
\subseteq M(S_V)$ and is a Hopf-*-homomorphism if and only if there
exists a non-degenerate representation $\psi$ of $(\hat{S}_U)_p$ on
${\cal L}(H)$ such that $(id\otimes \phi)(U'') = (\psi \otimes
id)(V')$. 
\par
\medskip
\noindent {\bf Lemma} 3.8{\bf :} There is only one $V$-$I_{{\bf
C}}$-birepresentation.
\par
\noindent {\bf Proof:} Let $\phi $ be a Hopf-*-homomorphism from 
$(S_{V})_{p}$ to ${\bf C}$. Then $\phi \in (S_{V})^{*}_{p}$ is an
idempotent. Suppose that $\chi $ is the homomorphism from 
$(S_{V})^{*}_{p}$ to 
$M(\hat{S}_{V})$ as defined in $[5,
A6]$. Then $\chi (\phi ) = (id\otimes \phi )(V^{\prime\prime})$ is both an 
idempotent and a unitary. Hence $\chi (\phi ) = 1.$
But since $\chi $ is unital and injective (see $[5, A6]$), $\phi $ is 
the co-identity of $(S_{V})_{p}$.
\par
\medskip
\noindent {\bf Definition} 3.9{\bf:} Let $W$ be a $C^*$-multiplicative 
unitary. Let $X$ be a $U$-$V$-birepresentation 
and $Y$ a $V$-$W$-birepresentation. 
Then the unitary $Z$ as given in Lemma 2.5 (if exists) is called 
the {\it composition} of $X$ and $Y$.
\par\medskip
\noindent {\bf Remark} 3.10{\bf:} (a) By Lemma 1.5, the composition $Z$
(if exists) is a $U$-$W$-birepresentation.
\par
\noindent (b) In the above setting, if $V$ is 
basic, then $Z= (\rho _{X}\otimes L_{Y})(V_{p})$ exists.
\par
\medskip
It is not clear for the moment how to relate $L_{Z}$ to $L_{X}$ and $L_{Y}$ 
i.e. how to
define ``composition'' of $L_{X}$ and $L_{Y}$. We will deal with this in 
the next section.
\par
\medskip
\noindent {\bf Example} 3.11{\bf :}(a) $X=V$ is a $V$-$V$-birepresentation.
\par
\noindent (b) $X=I_{H\otimes K}$ is a $U$-$V$-birepresentation. Hence the 
collection of $U$-$V$-birepresentations is non-empty. 
\par
\noindent (c) Let $G$ and $H$ be locally compact groups and $U=V_{H}$ and 
$V=V_{G}$ (where $V_{G}\xi (s,t) =
\xi (ts,t)$ for all $s,t\in G$ and $\xi \in L^{2}(G\times G))$. If 
$\phi $ is a group homomorphism from $G$ to
$H$, then $X\in L^{2}(H\times G)$ defined by $X\eta (r,s) = \eta (\phi 
(s)r,s)$ is a $U$-$V$-birepresentation
such that the map $L_{X}$ from $(S_{U})_{p} = C_{0}(H)$ to $M(S_{V}) = 
C_{b}(G)$ is the *-homomorphism
defined by $\phi $. In fact, for any $\xi ,\eta \in L^{2}(H)$, the map $g$ 
defined by $g(s) = (\omega _{\xi ,\eta }\otimes id)(U)(s) =
\int \xi (t)\overline{\eta (st)}dt$ is in $C_{0}(H)$ and $L_{X}(g)(r) = 
(\omega _{\xi ,\eta }\otimes id)(X)(r) = \int \xi (t)\overline{\eta 
(\phi (r)t)}dt = g(\phi (r))$. 
Note that by Theorem 3.6, $U$-$V$-birepresentations are 
precisely group homomorphisms in this case.
\par
\medskip
As a corollary of Proposition 1.9, we have the following:
\par
\medskip
\noindent {\bf Lemma} 3.12{\bf :} Let $\psi $ and $\phi $ be Hopf 
*-homomorphisms from $(S_{U})_{p}$ to $M(S_{V})$ and from
$(\hat{S}_{V})_{p}$ to $M(\hat{S}_{U})$ respectively. Let $\delta _{U}$ and 
$\hat{\delta }_{V}$ be the co-multiplications on $(S_{U})_{p}$
and $(\hat{S}_{V})_{p}$ respectively. Then $\epsilon = (id\otimes \psi )
\circ \delta _{U}$ is a coaction on $(S_{U})_{p}$ by $S_{V}$ and
$\hat{\epsilon } = (id\otimes \phi )\hat{\delta }_{V}$ is a coaction on 
$(\hat{S}_{V})_{p}$ by $\hat{S}_{U}$.
\par
\medskip
Let $X$ be a $U$-$V$-birepresentation. Then, by Lemma 3.12, $X$ induces 
coactions
$\epsilon _{X}$ and $\hat{\epsilon }_{X}$ on $(S_{U})_{p}$ and 
$(\hat{S}_{V})_{p}$ by $S_{V}$ and $\hat{S}_{U}$ respectively. Moreover, 
these coactions are "mutual" in the following sense:
\par
\medskip
\noindent {\bf Definition} 3.13{\bf :} Let $(S, \delta _{S})$ and $(T, 
\delta _{T})$ be two Hopf $C^{*}$-algebras. A coaction $\epsilon $
on $S$ by $T$ is said to be a {\it mutual coaction} if $(\delta _{S}\otimes 
id)\circ \epsilon = (id\otimes \epsilon )\circ \delta _{S}$.
\par
\medskip
We can now add one more equivalence to Theorem 3.6.
\par
\medskip
\noindent {\bf Proposition} 3.14{\bf :} Let $\epsilon $ be a mutual coaction 
of $S$ by $T$. If $S$ has a co-identity
$E$, then $\psi =(E\otimes id)\circ \epsilon $ is a Hopf *-homomorphism from 
$S$ to $M(T)$ such that $\epsilon =
(id\otimes \psi )\circ \delta _{S}$. Hence, $U$-$V$-birepresentations are in 
one to one correspondence with
mutual coactions of $(S_{U})_{p}$ by $S_{V}$ (and also with mutual coactions 
of $(\hat{S}_{V})_{p}$ by
$\hat{S}_{U}$).
\par
\noindent {\bf Proof:} First note that $\delta _{T}\circ \psi = 
(E\otimes id\otimes id)(\epsilon \otimes id)\epsilon = (\psi \otimes id)
\circ \epsilon $. Moreover, $(\epsilon \otimes \epsilon )\circ \delta _{S} = 
(\epsilon \otimes id\otimes id)(\delta _{S}\otimes id)\epsilon $. 
Hence, $(\psi \otimes \psi )\circ \delta _{S} = (E\otimes id\otimes 
id)(\epsilon \otimes id)(id\otimes E\otimes id)(\delta _{S}\otimes 
id)\epsilon =
(\psi \otimes id)\circ \epsilon $. Finally, $(id\otimes \psi )\circ \delta 
_{S} = (id\otimes E\otimes id)(\delta _{S}\otimes id)\epsilon = 
\epsilon $. Now by $[5, 2.1], (S_{U})_{p}$
has a co-identity and the second part follows from Theorem 3.6.
\par
\medskip
It is natural to ask what is the relation between the crossed product
of $\epsilon_X$ and that of $\hat \epsilon _X$. Before we compare
these two crossed products, let us first give the following lemmas.
\par
\medskip
\noindent {\bf Lemma} 3.15{\bf :} Let $B$ be a $C^{*}$-algebra and let 
$\phi $ and $\mu $ be *-homomorphisms from
$(S_{U})_{p}$ and $(\hat{S}_{V})_{p}$ respectively to $M(B)$. Then $(\phi , 
\mu )$ is a covariant pair for
$((S_{U})_{p}, S_{V}, \epsilon _{X})$ if and only if $((id\otimes \phi 
)(U^{\prime\prime}), (\mu \otimes id)(V^\prime ))$ is a covariant pair
for $X$ in the sense of Definition 1.6.
\par
\noindent {\bf Proof:} Let $u = (id\otimes \phi )(U^{\prime\prime})$ and 
$v = (\mu \otimes id)(V^\prime )$. Then, by definition, 
$(\phi , \mu )$ is a covariant pair
for $((S_{U})_{p},S_{V},\epsilon _{X})$ if and only if for any $\omega \in 
{\cal L}(K)_{*}$, {$(\phi \otimes id)(id\otimes L_{X})\delta 
_{U}[(\omega \otimes id)(U^{\prime\prime})]=$}
$(\mu \otimes id)(V^\prime )[\phi ((\omega \otimes 
id)(U^{\prime\prime}))\otimes 1](\mu \otimes 
id)(V^\prime )^{*}$. This is
 the case if and only if
$(\phi \otimes id)(\omega \otimes id\otimes id)(U_{12}^{\prime\prime}
X_{13})\cdot 
v = v\cdot (\omega \otimes id\otimes id)(u_{12})$. Now the left hand side 
equals $(\omega \otimes id\otimes id)(u_{12}X_{13}v_{23})$ 
while the right hand side 
is $(\omega \otimes id\otimes id)(v_{23}u_{12})$. Therefore,
the lemma follows from the fact that ${\cal L}(H)_{*}$ separates points of 
$S_{U}$.
\par
\medskip
\noindent {\bf Lemma} 3.16{\bf :} Let $B$, $\phi $ and $\mu $ be the same 
as in the previous lemma. Then $(\phi ,\mu )$
is a covariant pair for $((S_{U})_{p}, S_{V}, \epsilon _{X})$ if and only if 
$(\mu , \phi )$ is a covariant
pair for $((S_{V^{\top }})_{p}, S_{U^{\top }}, \epsilon _{X^{\top }})$.
\par
\noindent {\bf Proof:} Let $u=(id\otimes \phi )(U^{\prime\prime})$ and 
$v=(\mu \otimes id)(V^\prime )$ as in the previous lemma. Suppose
that $(\phi , \mu )$ is a covariant pair for $((S_{U})_{p}, S_{V}, \epsilon 
_{X})$. Then $u_{12}X_{13}v_{23}=v_{23}u_{12}$.
Now let $y = (id\otimes \mu )(V^{\top \prime\prime})$ and $z = (\phi 
\otimes id)(U^{\top \prime} )$. It is required to show that
$y_{12}X^{\top }_{13}z_{23}=z_{23}y_{12}$. In fact, $y=v^{*\sigma }$ 
and $z=u^{*\sigma }$ (where $\sigma $ is the flip of variables).
Thus, $y_{12}X^{\top }_{13}z_{23} = (\mu \otimes id)(V^{\prime 
*})_{21}X^{*}_{31}(id\otimes \phi )(U^{\prime\prime *})_{32} = 
(u_{32}X_{31}v_{21})^{*} = (v_{21}u_{32})^{*}$ (this 
is true by flipping the first and the third variables) and so 
$y_{12}X^{\top }_{13}z_{23} = z_{23}y_{12}$. 
The proof for the converse is the same.
\par\medskip
Actually, the crossed product of $\epsilon _X$ is the same as that of
the opposite of $\hat\epsilon_X$ i.e. $\epsilon_{X^\top}$. (Note that
$S_{U^\top} = (\hat{S})_U^{op}$).
\par
\medskip
\noindent {\bf Proposition} 3.17{\bf :} 
$(S_{U})_{p}\times _{\epsilon _{X},\max }\hat{S}_{V} \cong 
({S}_{V^{\top}})_{p}\times _{{\epsilon }_{X^{\top}},\max }
\hat{S}_{U^{\top}}$.
\par
\medskip
By this proposition and Corollary 2.12, we have:
\par
\medskip
\noindent {\bf Corollary} 3.18{\bf :} If $\delta $, 
$\delta ^{\top}$ and $\hat{\delta }$ are co-multiplications on
$(S_{V})_{p}$, $(S_{V^{\top}})_{p}$ and $(\hat{S}_{V})_{p}$ respectively, 
then $(S_{V})_{p}\times 
_{\delta 
,\max }(\hat{S}_{V})_{p} \cong (S_{V^{\top}})_{p}\times _{\delta ^{\top},
\max }(\hat{S}_{V^{\top}})_{p} \cong (\hat{S}_{V})_{p}
\times _{\hat{\delta}^{op},\max }(S_{V})_p^{op}$.
\par
\medskip
Now for any $U$-$V$-birepresentation $X$, we obtain a $
C^{*}$-algebra $C^{*}(X) =
(S_{U})_{p}\times _{\epsilon _{X},\max }\hat{S}_{V} = (\hat{S}_{V})_{p}
\times _{\hat{\epsilon}_{X^{\top}},\max }S_{U}^{op}$ which has coactions by 
$(\hat{S}_{V})_{p}$ 
and by $(S_{U})_{p}$ respectively (see Proposition 2.13) such
that the canonical maps $\mu $ and $\phi $ from $(\hat{S}_{V})_{p}$ and 
$(S_{U})_{p}$ respectively to $M(C^{*}(X))$ are
equivariant (see Proposition 2.13).
\par
\medskip
\noindent {\bf Remark} 3.19{\bf :} By Proposition 1.10, we 
obtain a map $\pi _{0}$ from 
$C^{*}(X)$ to $M(\hat{S}_{V}\times _{\delta _{V},\max }S_{V})$ and
hence a representation $\pi $ of $C^{*}(X)$ on H. Similarly, we have a 
representation $\tau$
of $C^{*}(X)$ on K. Moreover, $L_{X} = \pi \circ \phi $ and 
$\rho _{X} = \tau\circ \mu $. 
In fact, if $\psi $ is the canonical map from $S_{V}$ to
$M(S_{V}\times _{\delta _{V},r}\hat{S}_{V})$ (which equals ${\cal L}(H))$, 
then $\pi \circ \phi = \psi \circ L_{X} = L_{X}$ (since $\psi = L_{V})$.
\par\medskip
\par
\medskip
\noindent {\bf 4. Lifting of birepresentations}
\par
\medskip
In this section, we assume that all multiplicative unitaries are 
basic. We will show that any birepresentation 
$X\in M(\hat{S}_{U}\otimes S_{V})$ is the image of a unitary 
$(\hat{S}_{U})_{p}$-$(S_{V})_{p}$ birepresentation 
$X_{p}$ (see Definition 1.4). Consequently, we can lift 
any Hopf *-homomorphism $\phi :(S_{U})_{p}\rightarrow M(S_{V})$ to 
a Hopf *-homomorphism $\phi ^\prime 
:(S_{U})_{p}\rightarrow M(S_{V})_{p}$. First of all, let 
$X$ be a
$U$-$V$-birepresentation and define $X^\prime =(id\otimes L_{X})(U_{p})$ and 
$X^{\prime\prime}=(\rho _{X}\otimes id)(V_{p})$. Then:
\par
\medskip
\noindent {\bf Lemma} 4.1{\bf :} $X^\prime $ and $X^{\prime\prime}$ are 
unitary $(\hat{S}_{U})_{p}$-$S_{V}$-birepresentation and unitary
 $\hat{S}_{U}$-$(S_{V})_{p}$-birepresentation respectively such that $(\rho 
_{U}\otimes id)(X^\prime ) = X =
(id\otimes L_{V})(X^{\prime\prime})$.
\par
\noindent {\bf Proof:} Since $U_{12}^\prime (U_{p})_{13}U_{23}^
{\prime\prime} = U
_{23}^{\prime\prime}U_{12}^\prime $, we have $U_{12}^\prime 
X_{13}^\prime X_{23} = X_{23}U_{12}^\prime $. Thus, by
Lemma 1.5 and the fact $X$ is a representation of $V$, $X^\prime $ is a 
unitary $(\hat{S}_{U})_{p}$-$S_{V}$-birepresentation. The second part of 
the lemma is clear. 
\par
\medskip
\noindent {\bf Remark} 4.2{\bf :} (a) Note that by Proposition 2.11, 
$X^\prime $ and $X^{\prime\prime}$ are uniquely determined by $X$.
\par
\noindent(b) Since $X^\prime $ is a 
representation of $V$, we can define a map $\rho _{X^\prime }$
from $(\hat{S}_{V})_{p}$ to $M((\hat{S}_{U})_{p})$ such that $X^\prime = 
(\rho _{X^\prime }\otimes id)(V^\prime )$. Moreover, by a similar
argument as the proof of Proposition 3.2, $\rho _{X^\prime }$ is a 
Hopf *-homomorphism. Similarly, we 
have a Hopf *-homomorphism $L_{X^{\prime\prime}}$ from $(S_{U})_{p}$ to 
$M((S_{V})_{p})$ such that $X^{\prime\prime} =
(id\otimes L_{X^{\prime\prime}})(V^{\prime\prime})$.
\par
\noindent (c) Both $L_{X^{\prime\prime}}$ and $\rho _{X^\prime }$ preserve 
antipodes (by a similar argument as in the proof
of Proposition 3.2). Hence {\it any Hopf-*-homomorphism from 
$(S_{U})_{p}$ to $M((S_{V})_{p})$ preserves antipodes}. 
\par
\noindent (d) Both $L_{X^{\prime\prime}}$ and $\rho _{X^\prime }$ preserve 
co-identities in the following sense: if $E_{V}$ is
the co-identity of $(S_{V})_{p}$, then $E_{V}\circ L_{X^{\prime\prime}}$ is 
the co-identity of $(S_{U})_{p}$. This
follows directly from Lemma 3.8. Thus {\it any Hopf-*-homomorphism from 
$(S_{U})_{p}$ to $M((S_{V})_{p})$ preserves co-identities}. 
\par
\medskip
\noindent {\bf Lemma} 4.3{\bf :} Let $X$ and $Y$ be the unitaries 
as in Definition 3.9 
and let $Z=X\circ $Y (see Definition 3.9). 
Then $L_{Z} = L_{Y}\circ L_{X^{\prime\prime}}$
and $\rho _{Z} = \rho _{X}\circ \rho _{Y^\prime }$.
\par
\noindent {\bf Proof:} We first note that $L_{X^{\prime\prime}}((\omega 
\otimes id)(U^{\prime\prime})) = (\omega \otimes id)(X^{\prime\prime}) = 
(\omega \circ \rho _{X}\otimes id)(V_{p})$. Thus,
$L_{Y}\circ L_{X^{\prime\prime}}((\omega \otimes id)(U^{\prime\prime})) = 
(\omega \otimes id)(\rho _{X}\otimes L_{Y})(V_{p}) = (\omega \otimes 
$id)(Z). The proof for $\rho _{Z} = \rho _{X}\circ \rho _{Y^\prime }$
is the same.
\par
\medskip
\noindent {\bf Lemma} 4.4{\bf :} $(id\otimes 
L_{X^{\prime\prime}})(U_{p}) = 
(\rho _{X^\prime }\otimes id)(V_{p})$.
\par
\noindent {\bf Proof:} Let $X_{1}=(id\otimes L_{X^{\prime\prime}})(U_{p})$ 
and $X_{2}=(\rho _{X^\prime }\otimes id)(V_{p})$. First notice that both 
$X_{1}$ and
$X_{2}$ are unitary corepresentations of $(S_{V})_{p}$ (as 
$L_{X^{\prime\prime}}$ is a Hopf $^{*}$-homomorphism).
Moreover, $(id\otimes L_{V})(X_{1})= (id\otimes L_{X})(U_{p}) = X^\prime $ 
and $(id\otimes L_{V})(X_{2}) = (\rho _{X^\prime }\otimes id)(V^\prime )= 
X^\prime $.
Hence, the lemma follows from Proposition 2.11.
\par
\medskip
\noindent {\bf Definition} 4.5{\bf :} Let $X$ be a $U$-$V$-birepresentation. 
Then $X_{p} = (id\otimes L_{X^{\prime\prime}})(U_{p}) =
(\rho _{X^\prime }\otimes id)(V_{p})$ is called the {\it lifting} of $X$ in 
$M((\hat{S}_{U})_{p}\otimes (S_{V})_{p})$.
\par
\medskip
\noindent {\bf Remark} 4.6{\bf :} $X_{p}$ is a unitary $(\hat{S}
_{U})_{p}$-$(S_{V})_{p}$-birepresentation by Remark 2.4(b).
Moreover, $(id\otimes L_{V})(X_{p}) = X^\prime $ and $(\rho _{U}\otimes 
id)(X_{p}) = X^{\prime\prime}$.
\par
\medskip
\noindent {\bf Lemma} 4.7{\bf :} $L_{Z^{\prime\prime}} = 
L_{Y^{\prime\prime}}\circ 
L_{X^{\prime\prime}}\hbox{ and }\rho _{Z^\prime } = \rho _{X^\prime }\circ 
\rho _{Y^\prime }$.
\par
\noindent {\bf Proof:} For any $\omega \in {\cal L}(H)_{*}$, we have 
$L_{Y^{\prime\prime}}\circ L_{X^{\prime\prime}}[(\omega \otimes 
id)(U^{\prime\prime})] = (\omega \otimes 
L_{Y^{\prime\prime}})(X^{\prime\prime}) =
(\omega \otimes id)(\rho _{X}\otimes L_{Y^{\prime\prime}})(V_{p}) = (\omega 
\otimes id)(\rho _{X}\circ \rho _{Y^\prime }\otimes id)(W_{p}) = (\omega 
\otimes id)(\rho _{Z}\otimes id)(W_{p}) = (\omega \otimes 
id)(Z^{\prime\prime})$.
The proof for $\rho _{Z^\prime } = \rho _{X^\prime }\circ \rho _{Y^\prime }$ 
is the same.
\par
\medskip
\noindent {\bf Lemma} 4.8{\bf :} Let $X$ be a $U$-$V$-birepresentation 
and $Y$ a 
$V$-$W$-birepresentation. Let
$Z = X\circ $Y. Then there are equivariant maps from $C^{*}(Z)$ to 
$M(C^{*}(X))$ and to $M(C^{*}(Y))$ (see Remark 3.19 for the definition 
of $C^*$(X)).
\par
\noindent {\bf Proof:} We first show that $L_{X^{\prime\prime}}$ from 
$(S_{U})_{p}$ to $M[(S_{V})_{p}]$ is equivariant with
respect to the coactions $\epsilon _{Z}$ and $\epsilon _{Y}$ respectively.
 In fact, $\epsilon _{Y}\circ L_{X^{\prime\prime}} =
(id\otimes L_{Y})\circ \delta _{V}\circ L_{X^{\prime\prime}} = (id\otimes 
L_{Y})\circ (L_{X^{\prime\prime}}\otimes L_{X^{\prime\prime}})\circ \delta 
_{U} = (L_{X^{\prime\prime}}\otimes id)\circ (id\otimes L_{Z})\circ \delta 
_{U} = (L_{X^{\prime\prime}}\otimes id)\circ \epsilon _{Z}$.
Hence, we obtain a non-degenerate map $\Psi $ from $C^{*}(Z)$ to 
$M(C^{*}(Y))$ (by $[5, 3.9])$.
Now let $\phi , \mu , \phi ^\prime $ and $\mu ^\prime $ be the canonical 
maps from $(S_{U})_{p}, (\hat{S}_{W})_{p}, (S_{V})_{p}$ and 
$(\hat{S}_{W})_{p}$ to $C^{*}(Z)$ and $C^{*}(Y)$ respectively. Then we 
have $(\Psi \otimes id)\overline{\epsilon }_{Z}(\phi (s)\mu (t)) =
(\Psi \otimes id)[(\phi (s)\otimes 1)(\mu \otimes id)\hat{\delta }_{W}(t)] = 
(\phi 
^\prime (L_{X^{\prime\prime}}(s))\otimes 1)\cdot (\mu ^\prime \otimes 
id)\hat{\delta }
_{W}(t) = \overline{\epsilon }_{Y}(\Psi \otimes id)(\phi (s)\mu (t))$
for any $s\in (S_{U})_{p}$ and $t\in (\hat{S}_{W})_{p}$ (where 
$\overline{\epsilon }_{Y}$ and $\overline{\epsilon }_{Z}$ are the dual 
coactions as defined
in Proposition 2.13). The map from $C^{*}(Z)$ to $M(C^{*}(X))$ is defined 
similarly by
considering $C^{*}(Z) = (\hat{S}_{V})_{p}\times _{\hat{\epsilon }_{X},\max 
}S_{U}$.
\par
\medskip
We summarise the equivalences of $U$-$V$-birepresentations as follows:
\par
\medskip
\noindent {\bf Theorem} 4.9{\bf :} There are one to one correspondences 
between the collections of the
following objects:
\par
\noindent $(a)$ $U$-$V$-birepresentations;
\par
\noindent $(b)$ Hopf *-homomorphisms from $(S_{U})_{p}$ to $M(S_{V})$ 
(respectively, from $(\hat{S}_{V})_{p}$ to $M(\hat{S}_{U})$);
\par
\noindent $(c)$ mutual coactions on $(\hat{S}_{V})_{p}$ by $(\hat{S}_{U})$ 
(respectively, on $(S_{U})_{p}$ by $(S_{V})$);
\par
\noindent $(a^\prime )$ unitary 
$(\hat{S}_{U})_{p}$-$(S_{V})_{p}$-birepresentations;
\par
\noindent $(b^\prime )$ Hopf *-homomorphisms from $(S_{U})_{p}$ to 
$M[(S_{V})_{p}]$ (respectively, from $(\hat{S}_{V})_{p}$ to 
$M[(\hat{S}_{U})_{p}]$);
\par
\noindent $(c^\prime )$ mutual coactions on $(\hat{S}_{V})_{p}$ by 
$(\hat{S}_{U})_{p}$ (respectively, on $(S_{U})_{p}$ by $(S_{V})_{p}$).
\par
\noindent In this case, the Hopf *-homomorphisms in $(b)$ and 
$(b^\prime)$ preserve antipodes automatically. Moreover, 
Hopf-*-homomorphisms in $(b^\prime)$ preserve coidentity.
\par
\medskip
\noindent {\bf Remark} 4.10{\bf :} 
(a) It is also clear from the above results 
that unitary $\hat{S}_{U}$-$S_{V}$-, $(\hat{S}_{U})_{p}$-$S_{V}$-, 
$\hat{S}_{U}$-$(S_{V})_{p}$- and $(\hat{S}_{U})_{p}$-$(S_{V})_{p}$- 
birepresentations are all the same. 
\par
\noindent (b) Note that there may not be a one to one correspondence between 
the set of Hopf *-homomorphisms from $S_{U}$ to $M(S_{V})$ and the sets in 
Theorem 4.9. 
In fact, they are in one to one correspondence if and only if $U$ is 
co-amenable. Note that if the trivial birepresentation induces a 
Hopf-*-homomorphism from $S_U$ to $M(S_V)$, then $S_U$ is counital (by 
Lemma 3.8) which implies that $S_U=(S_U)_p$ (by [5, A4]). 
\par
\medskip
\par
\medskip
\noindent {\bf 5. The category of basic multiplicative unitaries}
\par
\medskip
Let ${\cal M}$ be the metagraph with the collection of all basic 
multiplicative
unitaries as its objects and birepresentations as arrows such that given a 
$U$-$V$-birepresentation $X$, we denote $dom(X)=V$ and $cod(X)=U$. 
Then by results in sections 3 and 4, we have the following:
\par
\medskip
\noindent {\bf Proposition} 5.1{\bf :} ${\cal M}$ is a category with 
null object $I_{{\bf C}}$. 
It contains the category of all locally compact groups as a full 
subcategory.
\par
\medskip
More generally, ${\cal M}$ also contains the category ${\cal M}_{ca}$ 
(respectively, ${\cal M}_a$) of co-amenable
(respectively, amenable) multiplicative unitaries in ${\cal M}$ 
as a full subcategory. Moreover, ${\cal M}_{ca}$ is a strict monoidal 
category (and so is ${\cal M}_{a}$) as shown by the following lemma.
\par
\medskip
\noindent {\bf Lemma} 5.2{\bf :} Let $U$ and $V$ be co-amenable 
(respectively, amenable) $C^*$-multiplicative unitaries. 
Then $W = U\otimes V$ is also a co-amenable (respectively, amenable)
$C^*$-multiplicative unitaries (hence is also basic by Remark 2.4). 
\par
\noindent {\bf Proof:} We first show that $W$ is a $C^*$-multiplicative 
unitary. In fact, by [5, 2.4 \& 2.5], any representation $X$ of $W$ is of 
the form $X= Y_{13}Z_{24}$ for some representations $Y$ and $Z$ of $U$ and 
$V$ respectively. Hence $\hat S_X$ is a $C^*$-algebra which equals 
$\hat S_Y\otimes \hat S_Z$. Hence condition (ii) of Definition 2.1 holds as 
well. Now we will show that $W$ is co-amenable. 
Let $E_{U}, E_{V}$ and $E_{W}$ be the co-identities 
of $(S_{U})_{p}, (S_{V})_{p}$ and $(S_{W})_{p}$
respectively. Since $U$ and $V$ are co-amenable, $(S_{U})_{p} = S_{U}$ and 
$(S_{V})_{p} = S_{V}$ are
nuclear (see $[5, 3.6])$. Let $E = E_{U}\otimes E_{V} \in 
[(S_{U})_{p}\otimes (S_{V})_{p}]^{*} = S^{*}_{W}$. Then $(E\circ L_{W})
(\omega \otimes \nu )
= (E_{U}\otimes E_{V})[L_{U}(\omega )\otimes L_{V}(\nu )] = (id\otimes 
\omega )(I_{K})\otimes (id\otimes \nu )(I_{H}) = 
E_{W}(\omega \otimes \nu )$. Now $L^{*}_{W}(S^{*}_{W})$ is a
right ideal of $(S_{W})_{p}^*$ (see $[5, A4]$) containing the identity 
and therefore $L^{*}_{W}(S^{*}_{W}) = (S_{W})_{p}^*$.
\par
\medskip
Since Woronowicz $C^{*}$-algebras will give multiplicative 
unitaries of compact type. We can roughly say that ${\cal M}_a$ 
contains all Woronowicz $C^{*}$-algebras if
we identify all Woronowicz $C^{*}$-algebras that give the same 
multiplicative unitaries. Now we turn to subobjects and quotients.
\par
\medskip
\noindent {\bf Definition} 5.3{\bf :} Let $U$ and $V$ be basic 
multiplicative unitaries.
\par
\noindent (a) $V$ is said to be a {\it sub-multiplicative unitary} of $U$ 
if there exists a Hopf-*-homomorphism  
$L_{X^{\prime\prime}}$ from $((S_{U})_{p})$ onto $(S_{V})_{p}$.
\par
\noindent (b) $U$ is said to be a {\it quotient} of $V$ if there exists a 
Hopf-*-homomorphism $\rho _{X^\prime }$ from $((\hat{S}_{V})_{p})$ onto 
$(\hat{S}_{U})_{p}$.
\par
\noindent (c) A $U$-$V$-birepresentation $X$ is said to be an 
{\it isomorphism} if there exist a
$V$-$U$-birepresentation $Y$ such that $X\circ Y = U$ and $Y\circ X$ =V.
\par
\medskip
\noindent {\bf Remark} 5.4{\bf :} Let $U$ and $V$ be basic 
multiplicative unitaries. Then an isomorphism between $U$ and $V$ 
is equivalent to the existence of two Hopf-*-homomorphisms $\psi $ and 
$\phi $ from 
$(S_{U})_{p}$ to $M((S_{V})_{p})$ and from $(S_{V})_{p}$ to $M((S_{U})_{p})$ 
respectively such that $\psi \circ \phi = id$ and $\phi \circ \psi =id$. 
Hence, by Lemma 1.3, isomorphisms between $U$ and $V$ are equivalent to 
Hopf-*-isomorphisms between $(S_{U})_{p}$ and $(S_{V})_{p}$. 
Moreover, if $S_{U} \cong S_{V}$ as Hopf 
$C^{*}$-algebras, then $U$ is isomorphic to $V$. 
\par
\medskip
\noindent {\bf Lemma} 5.5{\bf :} Let $U$, $V$ and $W$ be basic 
multiplicative unitaries. 
Let $X$ and $Y$
be a $U$-$V$-birepresentation and a $V$-$W$-representation respectively. 
Then $X\circ Y =I$ if and only if $\rho _{Y^\prime }$ map 
$(\hat{S}_{W})_{p}$ into the fixed point
algebra of $\hat{\epsilon}_{X}$ in $M((\hat{S}_{V})_{p})$.
\par
\noindent {\bf Proof:} Let $Z=X\circ $Y. Note that 
$\hat{\epsilon}_{X}\circ \rho _{Y^\prime } =
(id\otimes \rho _{X^\prime })\circ \hat{\delta}_{V}\circ \rho _{Y^\prime }$. 
If $Z=I$, then $\hat{\epsilon}_{X}\circ \rho _{Y^\prime } = 
(\rho _{Y^\prime }\otimes \rho _{Z^\prime })\circ \hat{\delta}_{W} = 
(\rho _{Y^\prime }\otimes
E\cdot 1)\circ \hat{\delta}_{W} = (\rho _{Y^\prime }\otimes 1)$ (where E is 
the coidentity of $(\hat{S}_{W})_{p}$). Conversely, suppose that
$\rho _{Y^\prime }(\hat{S}_{W})_{p}\subseteq M((\hat{S}_{V})_{p})
^{\hat{\epsilon}_{X}}$. Then for any $\omega \in {\cal L}(K)_{*}$ (where $K$ 
is the underlying Hilbert
space for $W$), $(id\otimes \rho _{X^\prime })\circ 
\hat{\delta}_{V}\circ \rho_{Y^\prime } ((id\otimes \omega )(W^\prime )) = 
\rho _{Y^\prime }((id\otimes
\omega )(W^\prime ))\otimes 1$. Now the right hand side of the equation
equals $(id\otimes \rho _{X^\prime })((id\otimes id\otimes \omega 
)(Y^\prime _{13}Y^\prime _{23})) = (id\otimes id\otimes \omega )(Y^\prime 
_{13}Z^\prime _{23})$ while the left
hand side is $(id\otimes id\otimes \omega )(Y^\prime _{13})$. Hence, 
$(id\otimes id\otimes \omega )(Y^\prime _{13}Z^\prime _{23}-Y^\prime _{13}) 
= 0$ for all
$\omega \in {\cal L}(K)_{*}$. Since ${\cal L}(K)_{*}$ separates points of 
${\cal L}(K)$, we have $Y^\prime _{13}Z^\prime _{23}=Y^\prime _{13}$ and 
thus $Z=I$ (as $Y^\prime $ is a unitary).
\par
\medskip
It is natural to ask whether we can define the kernel of a morphism. 
We don't know how to define it in general. 
However, by examining normal 
subgroups of discrete groups (see the Appendix) and suggested by the
above lemma, we try to define 
kernels of morphisms between basic multiplicative unitaries 
of discrete type as follows.
\par
\medskip
\noindent {\bf Definition} 5.6{\bf :} Let $U$, $V$ and $W$ be 
regular multiplicative unitaries of discrete type. Let $X$ be a 
$U$-$V$-birepresentation and $Y$ be a $V$-$W$-birepresentation. 
\par
\noindent (a) Then $Y$ is said to be a {\it kernel} of $X$ if 
$\rho _{Y^\prime }$ is an isomorphism from 
$(\hat{S}_{W})_{p}$ to the fixed
point algebra of $\hat{\epsilon}_{X}$ in $(\hat{S}_{V})_{p}$.
\par
\noindent (b) If $W$ is a submultiplicative unitary of $V$ through $Y$, 
then $W$ is said to be {\it normal} if $Y$ is a kernel of a morphism. 
\par
\medskip
\noindent {\bf Proposition} 5.7{\bf :} If the kernel of a morphism $X$ 
exists, then it 
is unique up to
isomorphism.
\par
\noindent {\bf Proof:} Let $(W_{1}, Y_{1})$ and $(W_{2}, Y_{2})$ be kernels 
of $X$. Then there exists a
Hopf-*-isomorphism between $(\hat{S}_{W_{1}})_{p}$ and 
$(\hat{S}_{W_{2}})_{p}$ (by definition). Now the proposition follows 
from Remark 5.4.
\par
\medskip
\noindent {\bf Proposition} 5.8{\bf :} Suppose that the kernel $(W, Y)$ 
of a $U$-$V$-birepresentation $X$ exists. 
Let $Y_{1}$ be a $V$-$W_{1}$-birepresentation 
such that $X\circ Y_{1} = I$. Then there exists a unique 
$W$-$W_{1}$-birepresentation 
$Z$ such that $Y_{1} = Y\circ Z$.
\par
\noindent {\bf Proof:} By Lemma 5.5, $\rho _{Y^\prime _{1}}$ induces a 
unique Hopf-*-homomorphism 
$\phi $ from $(\hat{S}_{W_{1}})_{p}$ to 
$(\hat{S}_{W})_{p}$ 
such that $\rho _{Y^\prime }\circ \phi = \rho _{Y^\prime _{1}}$ ($\phi $ 
is unique since $\rho _{Y^\prime }$ is injective). Now the 
proposition follows from Theorem 4.9. 
\par
\medskip
\noindent {\bf Remark} 5.9{\bf :} 
(a) Proposition 5.8 justifies the use of the term ``kernel''. 
\par
\noindent (b) The kernel of a morphism need not exist in general e.g. if 
$H$ is a closed subgroup of a compact group $G$, then the fixed point 
algebra $C(G)^{\alpha _{H}}$ (which equals $C(G/H)$) 
is not a Hopf $C^{*}$-subalgebra of $C(G)$ unless 
$H$ is normal (where ${\alpha _{H}}$ is the action of $H$ on 
$C(G)$ induced from the canonical action of $G$ on itself).
\par\medskip
\noindent {\bf Example} 5.10{\bf :} 
The only example about normal submultiplicative unitaries that 
we have, for the moment, is the following very simple one. Let $V$ be 
the product of $U$ and $W$, then $W$ is a normal submultiplicative 
unitary of $V$. 
\par
\medskip
\par\medskip
\noindent {\bf 6. An Imprimitivity Type Theorem for Multiplicative 
Unitaries of Discrete Type}
\par\medskip
Let $U$ be a regular multiplicative unitary of discrete type. 
$U$ is clearly co-amenable. 
Let $\phi_U$ be the Haar state on $\hat S_U$. 
If $\epsilon $ is a coaction on $A$ by $\hat S_{U}$ 
with fixed point algebra $A^{\epsilon}$, then 
$E = (id\otimes \phi_{U})\circ \epsilon$ is a conditional 
expectation from $A$ 
onto $A^{\epsilon}$. 
\par\medskip
In this section, we will give an imprimitivity type theorem for 
discrete type multiplicative unitaries. 
On our way to this, we found the following interesting fact from Lemma 
6.5: the set 
$\{ (\omega_{e,\xi}\otimes id)(U): \xi\in H \}$ generates $S_U$ if $U$ is 
of discrete type and $e$ is the co-fixed vector of $U$ (see [2, 1.8]). 

\par\medskip
Stimulated by [9, 2.2.16], 
we are going to use Watatani's $C^*$-basic 
construction (see [9, sections 2.1 and 2.2]) to prove the imprimitivity 
type theorem. 
We recall that if $A$ is a 
$C^*$-subalgebra of $B$ with a common unit and $E$ is faithful 
conditional expectation from $B$ to $A$, 
then the $C^*$-basic construction $C^*<B, e_A>$ is equal to ${\cal K(F)}$ 
where $\cal F$ is the completion of $B$ with respect to the norm defined by 
$E$ (see [9, 2.1.3 and 2.2.10]). Moreover, we recall the following result:
\par\medskip
\noindent {\bf Proposition} 6.1{\bf :} (Watatani) Let B be a unital 
$C^*$-algebra and $A$ is a subalgebra of $B$ that contains the unit of $B$. 
Let $E$ be a faithful conditional expectation 
from $B$ to $A$. If $B$ acts on a Hilbert space $H$ 
faithfully and $e$ is a projection on $H$ such that 

\noindent (i) $ebe = E(b)e$ for all $b\in B$ and 

\noindent (ii) the map that send $a\in A$ to $ae\in 
{\cal L}(H)$ is injective, 

\noindent then the norm closure of $BeB$ is isomorphic to 
$C^*<B,e_A>$ canonically.
\par\medskip
We now state the main theorem of this section.
Let $U$, $V$ and $W$ 
be multiplicative unitaries of discrete type such that $W$ is a normal 
sub-multiplicative unitary of $V$ with quotient $U$. Let $\epsilon'$ 
be the coaction on $(\hat S_V)_p$ by $(\hat S_U)_p$ as defined in 
Section 3. For technical reasons, we assume that $U$ is amenable. 
\par\medskip
\noindent {\bf Theorem} 6.2 {\bf :} $(\hat S_W)_p$ is strongly Morita 
equivalence to $(\hat S_V)_p\times_{\epsilon',r}S_U$.
\par\medskip
Note that there exist a faithful Haar state for $\hat{S}_U$ if $U$ is
of discrete type. 
Hence $U$ is biregular and irreducible (up to multiplicity). 
Since it is more convenient for us to work with the reduced crossed product
of the form $A\times_{\epsilon ,r}\hat{S}_U$, we will consider $\hat U$ 
instead of $U$. 
Let $\hat U$ be as defined in [2, 6.1]. By [2, 6.8], 
$\hat S_U \cong S_{\hat U}$ as Hopf $C^*$-algebras. Let $\epsilon$ be the 
coaction on $(\hat S_V)_p$ by $S_{\hat U}$ induced by $\epsilon'$ and 
let $\psi_U$ be the corresponding Haar state on $S_{\hat U}$. 
Then $(\hat S_V)_p\times_{\epsilon',r}S_U = 
(\hat S_V)_p\times_{\epsilon,r}\hat S_{\hat U}$. 
Let $E=(id\otimes\psi_U)\circ\epsilon$ 
be the conditional expectation 
from $(\hat S_V)_p$ to $(\hat S_W)_p$ as defined by the first 
paragraph of this section. Since 
$\psi_U$ is faithful, $E$ is faithful 
($\epsilon$ is injective since it is defined by a Hopf-*-homomorphism 
from $(\hat S_V)_p$ to $\hat S_U$ and $\hat S_U$ has a co-identity). We 
first give the following lemmas.
\par\medskip
\noindent {\bf Lemma} 6.3{\bf :} $(id\otimes id\otimes \psi_U)(\hat U_{12}
\hat U_{13}) = (id\otimes id\otimes \psi_U)(\hat U_{13})$. 
\par
\noindent {\bf Proof:} Since $\psi_U$ is the Haar state, $(id
\otimes \psi_U)\delta_{\hat U}(x) = \psi_U(x) \cdot 1$ 
for all $x\in S_{\hat U}$ (where 
$\delta_{\hat U}$ is the comultiplication on $S_{\hat U}$). 
Hence $(\omega \otimes id
\otimes \psi_U)(\hat U_{12}\hat U_{13}) = (\omega \otimes id
\otimes \psi_U)(\hat U_{13})$ for all $\omega\in {\cal L}(H_U)_*$ 
and the lemma follows immediately. 
\par\medskip
\noindent {\bf Lemma} 6.4{\bf :} Let $(\hat S_V)_p$ be 
faithfully represented on a 
Hilbert space $H$. Regard $\epsilon$ 
as an injective map from $(\hat S_V)_p$ to 
${\cal L}(H\otimes H_U)$ and let $e = 1\otimes p$ (where $p = (\phi_U
\otimes id)(U) = (id\otimes \psi_U)(\hat U) \in {\cal L}(H_U))$. 
Then $\epsilon$ and $e$ will satisfy the 
two conditions on Proposition 6.1.
\par
\noindent {\bf Proof:} Since for any $a\in (\hat S_W)_p$, 
$\epsilon (a) = a\otimes 1$, 
the map in (ii) of Proposition 6.1 will send $a$ to $a\otimes p$ and so 
is injective. 
We can formulate condition (i) in the following way: 
$(1\otimes p)\epsilon (b)
(1\otimes p) = (id\otimes \psi_U)\epsilon (b)\otimes p$ for 
any $b\in (\hat S_V)_p$. 
Now $(1\otimes p)\epsilon (b)(1\otimes p) = $ $(id\otimes id\otimes 
\psi_U \otimes \psi_U)(\hat U_{23}(\epsilon (b)\otimes 1\otimes 1)
\hat U_{24}) =$ $(id\otimes id\otimes \psi_U \otimes \psi_U)
(((id\otimes \delta_{\hat U})\epsilon (b)\otimes 1)\hat U_{23}\hat U_{34})$ 
(since $\epsilon$ satisfies the coaction identity). Thus, using 
Lemma 6.3, $(1\otimes p)\epsilon (b)(1\otimes p) =$ $(id\otimes id\otimes 
\psi_U \otimes \psi_U)(((id\otimes \delta_{\hat U})\epsilon (b)\otimes 1)
\hat U_{24})=$ $[(id\otimes id\otimes \psi_U)((id\otimes \delta_
{\hat U})\epsilon (b))](1\otimes p)=$
$(id\otimes\psi_U)(\epsilon (b))\otimes p$.
This proved the lemma. 
\par\medskip
\noindent {\bf Lemma} 6.5{\bf :} The set 
$P = \{ (\phi_U\cdot s \otimes id)(U): s\in \hat A_U \}$ is dense in $S_U$. 
Equivalently, $\{ (id \otimes \psi_U\cdot s)(\hat U): s\in A_{\hat U} \}$ 
is dense in $\hat S_{\hat U}$.
\par
\noindent {\bf Proof:} We first note that because $p=(\phi_U\otimes id)(U)$ 
is a minimum central projection, $p\in S_U$ ($p\cdot S_U = {\bf C}\cdot p$). 
Moreover, if $s=(id\otimes \omega)(U)$ then $(\phi_U\otimes id)((s\otimes 1)
U)= (\phi_U\otimes id)(id\otimes \omega\otimes id)(U_{12}U_{13}) = 
(\phi_U\otimes id)(id\otimes \omega\otimes id)((id\otimes \delta_U)U) = 
(\omega\otimes id)\delta_U(p)$. Note that $\delta_U(p)(x\otimes 1)\in 
S_U\otimes S_U$ (for any $x\in S_U$) and so 
$(\omega\otimes id)\delta_U(p)\in 
S_U$. Thus $P$ is a subset of $S_U$.
Let $t\in \hat S_U$ be such that $(\phi_U\cdot s)(t)=0$ for all $s\in 
\hat A_U$. Then $\phi_U(t^*t)=0$ (as $\hat A_U$ is dense in $\hat S_U$). 
Because $\phi_U$ is faithful, $P$ separates points of $\hat S_U$. 
Hence $P$ is $\sigma(\hat S_U^*, \hat S_U)$-dense 
in $\hat S_U^*$. Therefore, for any $f\in \hat S_U^*$, there exists a net 
${s_i}$ in $\hat A_U$ such that $\phi_U\cdot s_i$ converges to $f$ weakly. 
Note that $g(L_{U}(h)) = h(\rho_{U}(g))$ for all $g\in S_U^*$ and 
$h\in \hat S_U^*$ and that $\rho_{U}(S^*_U)$ 
is a dense subset of $\hat S_U$ (because $1\in \hat S_U$). 
Hence for any $\nu \in {\cal L}(H_U)_*$, there exists a net $a_i$ in 
$L_{U}(P)$ such that $g(a_i)$ converges to $g(L_{U}(\nu))$
for any $g\in S_U^*$. Therefore, the $\sigma (S_U, S_U^*)$-closure 
of $L_{U}(P)$ will contains $S_U$ and so $L_U(P)$ is norm dense in $S_U$ 
(because $L_U(P)$ is a convex subset, in fact a vector subspace, of $S_U$). 
\par\medskip
\noindent {\bf Lemma} 6.6{\bf :} 
Let the notation be the same as in Lemma 6.4. Then the linear span, $T$, of 
$\{ \epsilon (a)(1\otimes p)\epsilon (b): a,b\in (\hat S_V)_p\}$ is 
norm dense in $(\hat S_V)_p\times_{\epsilon,r}\hat S_{\hat U} = $ 
$(\hat S_V)_p\times_{\epsilon',r}S_U$. 
\par
\noindent {\bf Proof:} 
We first note that $T$ is a subset of 
$(\hat S_V)_p\times_{\epsilon,r}\hat S_{\hat U}$. 
Since $\epsilon$ is a coaction, 
$(1\otimes p) \epsilon (b) = 
(id\otimes id\otimes \psi_U)((\epsilon\otimes id)\epsilon(b)\hat U_{23})$. 
Therefore, $\epsilon (a)(1\otimes p)\epsilon (b) = 
(id\otimes id\otimes \psi_U)((\epsilon\otimes id)((a\otimes 1)\epsilon (b))
\hat U_{23})$. 
Now $((\hat S_V)_p\otimes 1)\epsilon (\hat S_V)_p =$ 
$(id\otimes \Phi)(((\hat S_V)_p\otimes 1)
\delta_V (\hat S_V)_p) =$ $(\hat S_V)_p\otimes S_{\hat U}$ 
(where $\delta_V$ is the comultiplication on $(\hat S_V)_p$ which is 
non-degenerate and $\Phi$ is the map from $(\hat S_V)_p$ to $S_{\hat U}$
that define $\epsilon$, $\Phi$ is surjective since $U$ 
is a quotient of $V$). 
Thus element of the form $(id\otimes id\otimes \psi_U)
((\epsilon \otimes id)(c\otimes s)\hat U_{23})$ ($c\in (\hat S_V)_p$ and 
$s\in S_{\hat U}$) can be approximated in norm by elements in $T$. Note that 
$(id\otimes id\otimes \psi_U)((\epsilon \otimes id)(c\otimes s)\hat U_{23})=$
$\epsilon (c)(1\otimes (id\otimes \psi_U\cdot s)(\hat U))$. Hence by Lemma 
6.5, $T$ is norm dense in $(\hat S_V)_p\times_{\epsilon,r}\hat S_{\hat U}$. 
\par\medskip
We can now prove the main theorem in this section very easily.
\par\medskip
\noindent {\bf Proof:} (Theorem 6.2) By Lemma 6.4 and Proposition 6.1 
(see the paragraph before Proposition 6.1 as well), $(\hat S_W)_p$ 
is strongly Morita equivalence to the closure of the linear span of the set 
$\{ \epsilon (a)(1\otimes p)\epsilon (b): a,b\in (\hat S_V)_p\}$ which, 
by Lemma 6.6, 
equals $(\hat S_V)_p\times_{\epsilon,r}S_U$.
\par\medskip
\noindent {\bf Remark} 6.7{\bf :}
It is believe that the amenability of $U$ can be removed.
\par
\medskip
\par
\medskip
\noindent {\bf Appendix}
\par
\medskip
The aim of this appendix is to give a $C^{*}$-algebraic characterisation of 
normal subgroups of discrete groups. Let $H$ be a discrete group and let 
$\Psi _{H}$ be the
canonical tracial state on $C^{*}(H)$. We first recall a well know fact 
about the fixed point algebra of a discrete coaction.
\par
\medskip
\noindent {\bf Lemma} A1{\bf :} Let $B$ be a $C^{*}$-algebra with a coaction 
$\epsilon $ by $C^{*}(H)$ and let $\Phi =
(id\otimes \Psi _{H})\circ \epsilon $. Then $\Phi (B)$ is the fixed point 
algebra $B^{\epsilon }$.
\par
\medskip
\noindent {\bf Theorem} A2{\bf :} Let $\varphi $ be a homomorphism from a 
discrete group G to a discrete group H. Let $N = Ker (\varphi )$ and 
$\epsilon _{H}$ be the coaction on $C^{*}(G)$ by $C^{*}(H)$ as given 
in Theorem 4.9. Then $C^{*}(N)$ is isomorphic to the fixed point algebra 
$C^{*}(G)^{\epsilon _{H}}$ of the coaction $\epsilon _{H}$.
\par
\noindent {\bf Proof:} By $[7, 4.1]$, the canonical map $j$ 
from $C^{*}(N)$ to 
$C^{*}(G)$ is injective. Therefore, we
need only to show that $j(C^{*}(N)) = C^{*}(G)^{\epsilon _{H}}$. It is clear 
that $j(C^{*}(N)) \subseteq 
C^{*}(G)^{\epsilon _{H}}$. Let $\Phi $ be the map as defined in Lemma A1 
with $B = C^{*}(G)$ and $\epsilon = \epsilon _{H}$. 
For any $t\in G$, $\Phi (u_{t}) = 
(id\otimes \Psi _{H})(u_{t}\otimes u_{\hat{t}})$ 
(where $\hat{t} = \varphi (t)$). 
Hence, $\Phi (u_{t})=0$ if $t\not\in N$ and $\Phi 
(u_{t})=u_{t}$ if $t\in N$. Now it is clear
that $C^{*}(G)^{\epsilon _{H}} \subseteq j(C^{*}(N))$ since 
$\bigoplus _{t\in G}{\bf C}\cdot u_{t}$ is a
dense subspace of $C^{*}(G)$.
\par
\medskip
\par
\medskip
\par
\medskip
\noindent {\bf REFERENCE:}
\par
\medskip
\noindent [1] S. Baaj, Representation reguliere du groupe quantique
des deplacements de Woronowicz, Asterisque (1995), no. 232, 11-48.
\par
\noindent [2] S. Baaj and G. Skandalis, Unitaires multiplicatifs et 
dualit\'e pour les
produits crois\'es de $C^{*}$-alg\`ebres, Ann. scient. \'Ec. Norm. Sup., 
$4^{e}$ s\'erie, t. 26
(1993), 425-488.
\par
\noindent [3] M. B. Landstad, J. Phillips, I. Raeburn and C. E. Sutherland, 
Representations of crossed products by coactions and principal bundles, 
Tran. of Amer. Math. Soc., vol 299, no. 2 (1987), 747-784.
\par
\noindent [4] V. de M. Iorio, Hopf $C^{*}$-algebras and locally compact 
groups, Pacific J.
Math. 87 (1980), 75-96.
\par
\noindent [5] C. K. Ng, Coactions and crossed products of Hopf 
$C^{*}$-algebras, Proc. London Math. Soc. (3) 72 (1996) 638-656.
\par
\noindent [6] I. Raeburn, On crossed products by coactions and their 
representation theory, Proc. London Math. Soc. (3), 64 (1992), 625-652.
\par
\noindent [7] M. A. Rieffel, Unitary representations of group extensions: 
An algebraic approach to the theory of Mackey and Blattner, in 
``Adv. in Math. Suppl. Stud.'', Vol 4, pp. 43-82, Academic Press, Orlando, 
FL, 1979.
\par
\noindent [8] S. Z. Wang, Free products of compact quantum groups, 
to appear in Comm. Math. Phys.
\par
\noindent [9] Y. Watatani, Index for $C^*$-algebras, Memoirs of Amer. Math. 
Soc. no. 242 (1990).
\par
\noindent [10] S. L. Woronowicz, Compact matrix pseudogroups, Comm. Math. 
Phys. 111(1987), 613-665.
\par
\noindent [11] S. L. Woronowicz, From multiplicative unitaries to 
quantum groups, Internat. J. Math. 7 (1996), no.1, 127-149.
\par
\medskip
\noindent MATHEMATICAL INSTITUTE, OXFORD UNIVERSITY, 24-29 ST. GILES, OXFORD 
OX1 3LB, UNITED KINGDOM.
\par
\noindent $e$-mail address: ng@maths.ox.ac.uk
\par
\end{document}